\documentclass[12pt]{article}
\usepackage{amsmath,amsthm,amssymb}
\usepackage{mathtools}
\newtheorem{theorem}{Theorem}[section] 
\newtheorem{lemma}{Lemma}[section]

\newtheorem{remark}{Remark}[section]

\numberwithin{equation}{section}

\newcommand{\D}{{\rm d}}
\newcommand{\dx}{\, \D x}

\newcommand{\ds}{\, \D s}
\newcommand{\dr}{\, \D r}
\newcommand{\dt}{\, \D t}

\newcommand{\dis}{\displaystyle}

\newcommand{\psp}{\,}

\newcommand{\rz}{\mathbb{R}}

\newcommand{\eps}{\varepsilon}

\newcommand{\klauf}{\left(\begin{array}}
\newcommand{\klzu}{\end{array}\right)}

\title{Liouville-type results in two dimensions for stationary points of functionals with linear growth}
\author{Michael Bildhauer \& Martin Fuchs}
\date{}

\newcommand{\reff}[1]{(\ref{#1})}

\usepackage{hyperref}

\hypersetup{
   pdfnewwindow=true,
   colorlinks=false,
   linkcolor=black,
   citecolor=green,
   filecolor=magenta,
   urlcolor=cyan,
}

\begin{document}

\parindent0em
\maketitle

\newcommand{\op}[1]{\operatorname{#1}}
\newcommand{\bv}{\op{BV}}
\newcommand{\mub}{\overline{\mu}}
\newcommand{\muhat}{\hat{\mu}}

\newcommand{\hypref}[2]{\hyperref[#2]{#1 \ref*{#2}}}
\newcommand{\hypreff}[1]{\hyperref[#1]{(\ref*{#1})}}

\newcommand{\ob}[1]{^{(#1)}}

\newcommand{\xh}{\Xi}
\newcommand{\oh}[1]{O\left(#1\right)}
\newcommand{\xn}{\underline{x}}
\newcommand{\yn}{\underline{y}}

\begin{abstract} 
We consider variational integrals of linear growth satisfying the condition of $\mu$-ellipticity
for some exponent $\mu >1$ and prove that stationary points $u$: $\rz^2 \to  \rz^N$ with the property
\[
\limsup_{|x|\to \infty} \frac{|u(x)|}{|x|} < \infty
\]
must be affine functions. 
The latter condition can be dropped in the scalar case together with appropriate
assumptions on the energy density providing an extension of Bernstein's theorem. 
\end{abstract}

\parindent0ex
%***************************************************************************************************
%***************************************************************************************************
\section{Introduction}\label{intro}
%***************************************************************************************************
%***************************************************************************************************
In this note we mainly present results of Liouville-type for entire solutions $u$: $\rz^2 \to \rz^N$ of the system
\begin{equation}\label{intro 1}
\op{div}\big[\nabla F(\nabla u)\big] = 0 \quad\mbox{on $\rz^2$}\, ,
\end{equation}
concentrating on the case of energy densities $F$: $\rz^{2N} \to \rz$ with linear growth.\\

To be precise we assume that $F$ is of class $C^2\big(\rz^{2N}\big)$ satisfying with constants $M$,
$\lambda$, $\Lambda >0$ and for some exponent $\mu >1$
\begin{eqnarray} 
\label{intro 2}
|\nabla F(Z)| &\leq & M \, ,\\
\lambda (1+|Z|)^{-\mu}|Y|^2 &\leq & D^2F(Z)(Y,Y) \leq \Lambda (1+|Z|)^{-1} |Y|^2
\label{intro 3}
\end{eqnarray}
for all $Y$, $Z\in \rz^{2N}$, where the first inequality in \reff{intro 3} expresses the fact that $F$ is a $\mu$-elliptic integrand. 
Note that \reff{intro 2} and \reff{intro 3} exactly correspond
to the requirements of Assumption 4.1 in \cite{Bi:1818} and as outlined in Remark 4.2 of this reference,
conditions \reff{intro 2} and \reff{intro 3} imply that $F$ is of linear growth in the sense that
\[
a |Z| - b \leq F(Z) \leq A |Z| + B \, ,\quad Z\in \rz^{2N}\, ,
\]
holds with constants $a$, $A>0$, $B$, $b$, $\geq 0$.\\

Note also that the ``minimal surface case'' is included by letting $F(Z) := (1+|Z|^2)^{1/2}$. In this case we have
the validity of \reff{intro 3} with the choice $\mu =3$, and two families of densities satisfying \reff{intro 2} and
\reff{intro 3} with prescribed exponent $\mu > 1$ are given by
\[
F(Z) := \left\{ \begin{array}{l}
\dis \int_0^{|Z|}\int_0^s (1+r)^{-\mu}\dr\ds\\[4ex]
\dis \int_0^{|Z|} \int_0^s (1+r^2)^{-\mu/2}\dr\ds
\end{array}\right\} \, , \quad Z \in \rz^{2N}\, .
\]
Our results on the behaviour of global solutions of the Euler equations \reff{intro 1} with $\mu$-elliptic
densities $F$ are as follows.\\

\begin{theorem}\label{intro theo 1}
Let $u \in C^2\big(\rz^2,\rz^N\big)$ denote a solution of \reff{intro 1} with density $F$ such that
\reff{intro 2} and \reff{intro 3} hold.
\begin{enumerate}
\item Suppose that in addition
\begin{equation}\label{intro 4}
\lim_{|x|\to \infty} \frac{|u(x)|}{|x|} = 0\, .
\end{equation}
Then $u$ is a constant function.
\item If the function $u$ has the property
\begin{equation}\label{intro 5}
\sup_{x\in \rz^2} |\nabla u(x)| < \infty \, ,
\end{equation}
then $u$ is affine.
\item If we have
\begin{equation}\label{intro 6}
\limsup_{|x|\to \infty} \frac{|u(x)|}{|x|} < \infty \, ,
\end{equation}
then the conclusion of $b)$ holds.\\
\end{enumerate}
\end{theorem}

\begin{remark}\label{intro remark 1}
\begin{enumerate}
\item Clearly \reff{intro 4} holds in the case that $u$ is a bounded solution, and evidently \reff{intro 5} implies \reff{intro 6}.
\item We do not know if there are versions of Theorem \ref{intro theo 1} for entire solutions
$u$: $\rz^n \to \rz^N$ of \reff{intro 1} in the case $n \geq 3$.
\item Our discussion of smooth solutions of the system \reff{intro 1} includes the
vector case $N>1$ for densities $F$ of linear growth. The existence of smooth solutions is known provided that
$\mu$ is not too large and provided that $F(Z) =f(|Z|)$.
 It is a challenging question whether the smoothness of solutions remains true (to some extend) if the second hypothesis
is dropped.\\
\end{enumerate}
\end{remark}

Before presenting the proof of Theorem \ref{intro theo 1} we wish to mention that there exists a variety of
Liouville-type theorems for entire solutions $u$: $\rz^n\to \rz^N$, $n \geq 2$, $N\geq1$, of systems of the form \reff{intro 1}
(and even for nonhomogeneous systems not generated by a density $F$) assuming that $F$ is of superlinear growth.
The interested reader should consult the references on this topic quoted for example in the
textbooks \cite{Gia:1983_1}, \cite{GM:2012_1}, \cite{GT:1998_1}, \cite{Gi:2003_1}, \cite{Jo:2013_1} and \cite{Ne:1983_1}.
A nice survey is also presented in \cite{Fa:2007_1}.\\

Besides this more general discussion the validity of Liouville theorems for harmonic maps between
Riemannian manifolds turned out to be a useful tool for the analysis of the geometric properties
of the underlying manifolds. Without being complete we refer to
\cite{Ch:1980_1}, \cite{EL:1988_1}, \cite{Hi:1982_1}, \cite{HJW:1980_1},  
\cite{HK:1972_1},  \cite{HKW:1977_1}, \cite{Sc:1984_1} and \cite{SY:1976_1}.\\

Liouville theorems are also of interest in the setting of fluid mechanics, where in the stationary case \reff{intro 1}
is replaced by a nonlinear variant of the Navier-Stokes equation with dissipative potential $F$ of superlinear
growth and the incompressibility condition $\op{div} u =0$ for the velocity field $u$: $\rz^n \to\rz^n$
has to be added. The validity of Liouville theorems has been established in the $2$-D-case, i.e.~for $n=2$, for instance
in the papers 
\cite{BFZ:2013_1}, \cite{Fu:2012_1}, \cite{Fu:2014_1}, \cite{FM:2019_1}, \cite{FZ:2012_1},
\cite{GW:1978_1}, \cite{JK:2014_1}, \cite{KNSS:2009_1}, \cite{Zh:2013_1} and \cite{Zh:2015_1}.
We like to mention that the case of potentials $F$ satisfying \reff{intro 2} and \reff{intro 3} is treated in
\cite{FM:2019_1} assuming $\mu < 2$.\\

As it stands, the conclusions of Theorem \ref{intro theo 1} b) and c) are in the spirit of Bernstein's theorem (see \cite{Be:1927_1})
for nonparametric minimal surfaces,
where in this particular setting conditions like \reff{intro 5} or \reff{intro 6} are seen to be superfluous.  
For completeness we specialize the Bernstein result obtained by Farina, Sciunzi and Valdinoci in Theorem 1.4 of their
paper \cite{FSV:2008_1} to the case of linear growth integrands.

\begin{theorem}\label{intro theo 2}
Consider a function $g \in C^2\big([0,\infty)\big)$ such that with constants $a_1$, $a_3$, $a_5 > 0$, $a_2$, $a_4 \geq 0$ we have for some exponent
$\mu \geq 1$
\begin{eqnarray}\label{intro 7}
& \dis g'(0) = 0 \, ,\quad g''(t) >0 \quad \mbox{for}\quad t > 0 \, ,\\[2ex]
& \dis a_1 t - a_2 \leq g(t) \leq a_3 t + a_4 \quad\mbox{for}\quad t \geq 0 \, ,& \label{intro 8}\\[2ex]
& \dis g''(t) \leq a_5 (1+t)^{-\mu} \quad\mbox{for}\quad t\geq 0 \, . \label{intro 9}
\end{eqnarray}
Let $F$: $\rz^2 \to \rz$, $F(Z) := g\big(|Z|\big)$, and consider a solution $u$: $\rz^2 \to \rz$ of \reff{intro 1} being of class $C^2$.
Then $u$ is an affine function provided that $\mu \geq 3$.
\end{theorem}

\begin{remark}\label{intro rem 2}
\begin{enumerate}
\item Note that the minimal surface case is included with the choices $g(t) = \sqrt{1+t^2}$ and $\mu =3$, moreover, we can cover the
examples stated in front of Theorem \ref{intro 1} provided that $\mu \geq 3$.

\item To our knowledge it is an unsolved problem, if Theorem \ref{intro theo 2} remains true for exponents $\mu \in (1,3)$.

\item Roughly speaking it follows from the work \cite{NN:1959_1} of J.C.C.~and J.~Nitsche that the Bernstein property fails for the equation
\[
0 = \op{div} \Bigg[\frac{g'\big(|\nabla u|\big)}{|\nabla u|} \nabla u\Bigg]\, ,
\]
if the density of $g$ is elliptic and of superlinear growth including even the nearly linear case $g(t) = t \ln(1+t)$, i.e.~there 
exist non-affine solutions $u$: $\rz^2 \to \rz$. However, the Nitsche criterion does not apply to integrands of linear growth
as considered in Theorem \ref{intro theo 2} (see Remark \ref{bern rem 1}).
\item From the identity
\begin{eqnarray*}
D^2F(Z)(X,X) &=& \frac{1}{|Z|} g'\big(|Z|\big)\Big[ |X|^2 - \frac{1}{|Z|^2} (X\cdot Z)^2\Big]\\ 
&&+ g''\big(|Z|\big) \frac{1}{|Z|^2} (X\cdot Z)^2\, , \quad X \,\, Z \in \rz^2\, ,\\
\end{eqnarray*}
it follows that (observing the boundedness of $g'$)
\[
\min \Bigg\{g''\big(|Z|\big), \frac{g'\big(|Z|\big)}{|Z|}\Bigg\} |X|^2 
\leq D^2F(Z)(X,X) \leq \Lambda \big(1+|Z|\big)^{-1} |X|^2 \, ,
\]
i.e.~the second inequality in \reff{intro 3} holds with some constant $\Lambda >0$. For $t\geq 1$ we have the lower bound
$g'(t)/t \geq c/t$, which by \reff{intro 9} means that in fact $g''\big(|Z|\big)$ measures the degree of ellipticity of $D^2F(Z)$. 
This shows that the integrand $F$ in general is not $\mu$-elliptic in the sense
of the first inequality from \reff{intro 3}: according to \reff{intro 9} the power $t^{-\mu}$ just acts as an upper bound for the values
$g''(t)$. Thus we have the ``Bernstein property'' for any density $F(Z) = g\big(|Z|\big)$ of linear growth and for which
$g''(t) = O(t^{-3})$ as $t\to \infty$.
\end{enumerate}
\end{remark}
%***************************************************************************************************
%***************************************************************************************************
\section{Proof of Theorem \ref{intro theo 1}, Part a)}\label{Pr}
%***************************************************************************************************
%***************************************************************************************************

In the weak formulation of \reff{intro 1}, i.e.~in the equation
\begin{equation}\label{Pr 1}
\int_{\rz^2} \nabla F(\nabla u): \nabla \varphi \dx = 0 \, , \quad \varphi \in C^1_{0}\big(\rz^2,\rz^N\big) \, ,
\end{equation} 
the function $\varphi$ is replaced by $\partial_\alpha \varphi$ ($\alpha \in \{1,2\}$ fixed), where now
$\varphi \in C^2_0\big(\rz^2,\rz^N)$ is assumed. With an integration by parts we obtain from \reff{Pr 1}
\begin{equation}\label{Pr 2}
\int_{rz^2} D^2F(\nabla u)\big(\partial_\alpha \nabla u, \nabla \varphi\big) \dx = 0 \, .
\end{equation} 
Now we choose $\varphi = \eta^2 \partial_\alpha u \in C^1_0\big(\rz^2,\rz^N\big)$ in \reff{Pr 2},
where $\eta \in C^1_0\big(\rz^2\big)$, $\op{spt}\eta \subset B_{2R}(0)$,
$\eta \equiv 1$ on $B_R(0)$, $0 \leq \eta \leq 1$, $|\nabla \eta| \leq c R^{-1}$. Then by
Cauchy-Schwarz's and Young's inequality we have (summation w.r.t.~$\alpha=1$, $2$)
\begin{eqnarray}\label{Pr 3}
\lefteqn{\int_{B_{2R}(0)} \eta^2 D^2F(\nabla u) \big(\partial_\alpha \nabla u,\partial_\alpha \nabla u\big)\dx}\nonumber\\
&\leq& c \int_{B_{2R}(0)} D^2F(\nabla u)\big(\nabla \eta \otimes \partial_\alpha u,\nabla \eta\otimes \partial_\alpha u\big)\dx \, . 
\end{eqnarray}
The hypotheses \reff{intro 2} and \reff{intro 3} yield
\begin{eqnarray}\label{Pr 4}
\int_{B_R(0)} \big(1+|\nabla u|\big)^{-\mu} |\nabla^2 u|^2 \dx &\leq&
c R^{-2} \int_{B_{2R}(0)-B_R(0)} \frac{|\nabla u|^2}{\sqrt{1+|\nabla u|^2}} \dx\nonumber\\
&\leq & cR^{-2} \int_{B_{2R}(0)-B_R(0)} |\nabla u| \dx  
\end{eqnarray}
and using the auxiliary inequality \reff{aux 1} of Lemma \ref{lem aux} given below we obtain for any $\eps >0$
\begin{eqnarray}\label{Pr 5}
\lefteqn{\int_{B_R(0)} \big(1+|\nabla u|\big)^{-\mu} |\nabla ^2 u|^2 \dx}\nonumber\\ 
&\leq& \frac{c}{R^2} \int_{B_{2R}(0)-B_R(0)} \Big[\eps + c(\eps)\big(\nabla F(\nabla u) - \nabla F(0)\big):\nabla u\Big]\dx \, .
\end{eqnarray}

With \reff{Pr 1} we also have
\begin{equation}\label{Pr 6}
\int_{\rz^2} \big(\nabla F(\nabla u)-\nabla F(0)\big): \nabla \varphi \dx = 0 \, , \quad \varphi \in C^1_{0}\big(\rz^2,\rz^N\big) \, ,
\end{equation} 
where we now choose $\varphi = \tilde{\eta}^2u$, $\tilde{\eta}\in C^1_0\big(\rz^2\big)$,
$\tilde{\eta} \equiv 1$ on $B_{2R}(0) - B_R(0)$, $\op{spt}\eta \subset B_{5R/2}(0)-\overline{B}_{R/2}(0)$, 
$0\leq \tilde{\eta}\leq 1$, $|\nabla \tilde{\eta}| \leq c/R$.\\

With this choice \reff{Pr 6} gives
\begin{eqnarray}\label{Pr 7}
\lefteqn{\int_{\rz^2} \big(\nabla F(\nabla u) - \nabla F(0)\big): \nabla u \tilde{\eta}^2 \dx}\nonumber\\
&=& - 2 \int_{\rz^2}\tilde{\eta} \big(\nabla F(\nabla u) - \nabla F(0)\big): (\nabla \tilde{\eta}\otimes u)\dx\nonumber\\
&\leq & c R^{-1} \int_{B_{5R/2}(0)-B_{R/2}(0)} |u| \dx \nonumber \\
& \leq & c R \sup_{B_{5R/2}(0)-\overline{B}_{R/2}(0)} |u| \, ,
\end{eqnarray}
where our assumption \reff{intro 2} is used.\\

By the definition of $\tilde{\eta}$ we obtain using \reff{Pr 7}
\begin{eqnarray}\label{Pr 8}
\lefteqn{ \int_{B_{2R}(0)-B_R(0)} \big(\nabla F(\nabla u) - \nabla F(0)\big):\nabla u\dx}\nonumber\\
& \leq &  \int_{\rz^2} \big(\nabla F(\nabla u) - \nabla F(0)\big): \nabla u \tilde{\eta}^2 \dx\nonumber\\
& \leq &  c  R \sup_{B_{5R/2}(0)-\overline{B}_{R/2}(0)} |u| \, .
\end{eqnarray}

If we insert \reff{Pr 8} into inequality \reff{Pr 5} and pass to the limit $R\to \infty$ recalling \reff{intro 4},
we obtain for any $\eps >0$
\[
\int_{\rz^2} \big(1+|\nabla u|\big)^{-\mu} |\nabla^2 u|^2 \dx \leq c \eps \psp ,
\]
hence $\nabla^2 u \equiv 0$ and therefore we find $A\in \rz^{2N}$, $a\in \rz^N$ such that
\[
u(x) = A x +a \, .
\]
Again we apply of the growth condition \reff{intro 4} and obtain $A=0$, hence the first part of 
Theorem \ref{intro theo 1} is established.\\

During the proof we made use of the elementary lemma

\begin{lemma}\label{lem aux}
Let $F\in C^2\big(\rz^{2N}\big)$ just satisfy the first inequality of \reff{intro 3} and let 
\[
\theta(r) := \frac{\lambda}{\mu -1} \big[ 1-(1+r)^{1-\mu}\big]\,  ,\quad r\geq 0\, .
\]
Then it holds for any $\eps >0$ and all
$Z\in \rz^{nN}$
\begin{equation}\label{aux  1}
|Z| \leq  \eps + \theta^{-1}(\eps) \big[\nabla F(Z) - \nabla F(0)\big] : Z \, .
\end{equation}
\end{lemma}

\mbox{}\\
\emph{Proof of Lemma \ref{lem aux}.} We fix $\eps > 0$.  If $|Z| \geq \eps$ then
\[
|Z| \theta \big(|Z|\big) \geq |Z| \theta (\eps) \psp ,
\]
which implies
\[
|Z| \leq \theta^{-1}(\eps) \, |Z|\,  \theta\big(|Z|\big)\, ,
\]
and if $Z\in \rz^{2N}$ is arbitrarily given, we have
\[
|Z| \leq \eps + \theta^{-1}(\eps)\, |Z| \, \theta\big(|Z|\big) \, .
\]
Moreover, 
\begin{equation}\label{aux 2}
\theta\big(|Z|\big)\,  |Z| \leq \big[ \nabla F(Z) - \nabla F(0)\big]:Z \
\end{equation}
easily follows from the first inequality in \reff{intro 3} as outlined in \cite{Bi:1818}, formula (1), p.~98.,
and \reff{aux 2} gives \reff{aux 1}.  \hspace*{\fill}$\Box$\\

\begin{remark}\label{P2 rem 1}
Clearly Lemma \ref{lem aux} is not limited to the case $n=2$ and without condition \reff{intro 3}
it would be sufficient to assume \reff{aux 2} for an
increasing non-negative function $\theta$: $[0,\infty) \to \rz$.
\end{remark}

%***************************************************************************************************
%***************************************************************************************************
\section{Proof of Theorem \ref{intro theo 1}, Parts b) and c)}\label{P2}
%***************************************************************************************************
%***************************************************************************************************

For Part b) we remark, that the idea of applying a Liouville argument to the derivatives of solutions, 
which are seen to solve an appropriate elliptic equation, has been successfully used by Moser \cite{Mo:1961_1}, Theorem 6,
with the result that entire solutions of the minimal surface equation
with bounded gradients in fact must be affine functions in any dimension $n \geq 2$.\\

In our setting, i.e.~for $n=2$ together with $N\geq 1$, one may just follow the arguments presented in \cite{Gia:1983_1},
Chapter III, p.~82, for an elementary proof essentially based on the ``hole-filling'' technique.\\

In Theorem \ref{intro theo 1}, Part b) turns out to be a corollary of Part c), which we now prove following some ideas given
in \cite{FZ:2012_1}.\\

As in the proof of the first part of Theorem \ref{intro theo 1} we obtain from \reff{Pr 3} the following variant of inequality \reff{Pr 4}
\begin{equation}\label{P2 1}
\int_{B_R(0)} D^2 F(\nabla u)\big(\partial_\alpha \nabla u,\partial_\alpha \nabla u\big) \dx
\leq cR^{-2} \int_{B_{2R}(0)-B_R(0)} |\nabla u| \dx 
\end{equation}
and, as outlined after \reff{Pr 4},  \reff{P2 1} gives for all $R >0$ and with the choice $\eps =1$
\begin{equation}\label{P2 2}
\int_{B_R(0)} D^2F(\nabla u) \big(\partial_\alpha \nabla u,\partial_\alpha \nabla u\big) \dx \leq
c \Big[ 1+ R^{-1} \sup_{B_{5R/2}(0)-B_{R/2}(0)} |u|\Big] \, .
\end{equation}
Inequality \reff{P2 2} shows, using \reff{intro 6}, 
\begin{equation}\label{P2 3}
\int_{\rz^2} D^2 F(\nabla u) \big(\partial_\alpha \nabla u, \partial_\alpha \nabla u\big) \dx < \infty \, .
\end{equation}
We finally claim that
\begin{equation}\label{P2 4}
\int_{\rz^2} D^2F(\nabla u) \big(\partial_\alpha \nabla u,\partial_\alpha \nabla u\big) \dx = 0 \, ,
\end{equation}
which gives $|\nabla^2 u| = 0$, hence the proof will be complete.\\

To prove \reff{P2 4} we again consider \reff{Pr 2} and choose $\varphi$ as done after this inequality. 
We obtain with $T_R := B_{2R}(0)-\overline{B}_{R/2}(0)$
using the Cauchy-Schwarz inequality
\begin{eqnarray*}
\lefteqn{\int_{\rz^2} D^2F(\nabla u) \big(\partial_\alpha \nabla u,\partial_\alpha \nabla u\big) \eta^2 \dx}\nonumber\\
&=&-2 \int_{T_R} D^2 F(\nabla u) \big(\eta \partial_\alpha \nabla u,\nabla \eta \otimes \partial_\alpha u\big)\dx\nonumber\\
&\leq & c \Bigg[\int_{T_R} \eta^2 D^2F(\nabla u) \big(\partial_\alpha \nabla u,\partial_\alpha \nabla u\big)\dx\Bigg]^{\frac{1}{2}}\nonumber\\
&& \cdot \Bigg[\int_{T_R} D^2 F(\nabla u) \big(\nabla \eta\otimes \partial_\alpha u,\nabla \eta \otimes \partial_\alpha u\big)\dx\Bigg]^{\frac{1}{2}}
\nonumber\\[2ex]
&=:& I_1(R) \cdot I_2(R)\, .
\end{eqnarray*}
We recall \reff{P2 3} which gives
\[
I_1(R) \to 0 \quad\mbox{as $R \to \infty$.}
\]
Assumption \reff{intro 3} yields the estimate
\[
I_2(R) \leq  c \Bigg[ R^{-2} \int_{T_R} |\nabla u|\dx \Bigg]^{\frac{1}{2}}\, .
\]
thus we obtain \reff{P2 4}, if we can prove 
\begin{equation}\label{P2 6}
\int_{B_R(0)} |\nabla u|\dx \leq c\big(1+R^2\big)\, .
\end{equation}
For \reff{P2 6} we use \reff{aux 1} (recall $\eta \equiv 1$ on $B_R(0)$) with the choice $\eps =1$, hence (compare
the derivation of \reff{Pr 7})
\begin{eqnarray*}
\int_{B_R(0)} |\nabla u|\dx &\leq & |B_R(0)| + 
c \int_{B_R(0)} \big[ \nabla F(\nabla u) - \nabla F(0)\big]: \nabla u \dx\\
&\leq &  |B_R(0)| + 
c \int_{B_2R(0)} \eta^2 \big[ \nabla F(\nabla u) - \nabla F(0)\big]: \nabla u \dx\\
&\leq &c \Big[R^2 + R \sup_{T_R} |u|\Big]\\
&=& c R^2 \Bigg[1 + \frac{1}{R} \sup_{T_R} |u| \Bigg]\, ,
\end{eqnarray*}
and our hypothesis \reff{intro 6} gives \reff{P2 4}, hence the proof of Theorem \ref{intro theo 1} is complete.
\hspace*{\fill}$\Box$\\

%***************************************************************************************************
%***************************************************************************************************
\section{Proof of Theorem \ref{intro theo 2}}\label{Bern}
%***************************************************************************************************
%***************************************************************************************************

We follow the lines of \cite{FSV:2008_1} by checking the hypotheses of Theorem 1.4 in this reference. We let
\begin{equation}\label{bern 1}
a(t) := \frac{g'(t)}{t}\, , \quad t > 0 \, ,
\end{equation}
and observe that on account of \reff{intro 7} the function $a$ continuously extends to $t=0$ by letting
$a(0) = g''(0)$. Obviously $a$ satisfies \reff{intro 2} from \cite{FSV:2008_1} ($g'$ is strictly increasing and thereby
positive on $(0,\infty)$ due to \reff{intro 7}), and since for any $t >0$ it holds (compare \reff{bern 1}) 
\begin{equation}\label{bern 2}
\lambda_1(t):= a(t) + t a'(t) = g''(t)
\end{equation}
we get (1.3) in \cite{FSV:2008_1}. At the same time assumption (A2) from \cite{FSV:2008_1} is obvious by formula \reff{bern 1}
and our requirements concerning $g$. Moreover, the stability condition (see (1.11) in \cite{FSV:2008_1}) follows from
\begin{equation}\label{bern 3}
\int_{\rz^2} \big( A(\nabla u) \nabla \phi\big) \cdot \nabla \phi \dx \geq 0
\end{equation}
for any $\phi \in C^1_0(\rz^2)$ with matrix ($Z \in \rz^2-\{0\}$)
\[
A_{ij}(Z) := |Z|^{-1} a'\big(|Z|\big) Z_i Z_j + a \big(|Z|\big) \delta_{ij}
\]
by observing that (recall \reff{bern 1}, \reff{intro 7})
\begin{eqnarray*}
\big(A(Y)Z\big) \cdot Z &=& |Y|^{-1} \Big[g''\big(|Y|\big) |Y|^{-1} - g'\big(|Y|\big) |Y|^{-2}\Big](Y\cdot Z)^2\\
&& + |Y|^{-1} g'\big(|Y|\big) |Z|^2\\
&\geq & g'\big(|Y|\big) \Big[|Y|^{-1} |Z|^2 - |Y|^{-3} (Y\cdot Z)^2\Big] \geq 0
\end{eqnarray*}
on account of $g'\big(|Y|\big) > 0$ for $Y\not= 0$.\\

It remains to check (1.17) and (1.18) in \cite{FSV:2008_1}: since $g$ is convex (see \reff{intro 7}) and of linear growth (compare \reff{intro 8})
the boundedness of $g'$ follows, hence we get (1.17) and by monotonicity $g'_\infty := \lim_{t\to \infty} g'(t)$
exists in $(0,\infty)$. By \reff{bern 2} the function $\lambda_1$ defined in \reff{bern 2} (see (1.5) of \cite{FSV:2008_1}) is just $g''$
so that (1.18) is a consequence of (1.9) and the aforementioned limit behaviour of $g'$ provided we assume $\mu \geq 3$.
From Theorem 1.4 in \cite{FSV:2008_1} it follows that $u(x) = \tilde{u}(\omega \cdot x)$ for some $\omega \in \rz^2$, $|\omega| =1$, 
and a function $\tilde{u}$: $\rz \to \rz$. Assuming $\omega =(1,0)$ (w.l.o.g.) we see that \reff{intro 1} implies
\[
\frac{\D}{\dt} \Bigg[\frac{1}{|\tilde{u}'|} g'\big(|\tilde{u}'|\big) \tilde{u}'\Bigg] = 0 \psp ,
\]
hence
\begin{equation}\label{bern 4}
|\tilde{u}'|^{-1} g'\big(|\tilde{u}'|\big)\tilde{u}' \equiv c
\end{equation}
for some $c \in \rz$.\\

\emph{Case 1, $c=0$.} Recalling $g'>0$ on $(0,\infty)$ equation \reff{bern 4} yields $\tilde{u}'\equiv 0$ and we are done.\\

\emph{Case 2, $c\not= 0$.} Then \reff{bern 4} shows $\tilde{u}'(t) \not= 0$ for any $t\in \rz$, thus
$g'\big(|\tilde{u}'|\big) \equiv |c|$ and in conclusion
\[
\tilde{u}'(t) \in \Big\{ -(g')^{-1}\big(|c|\big), (g')^{-1}\big(|c|\big) \Big\}
\]
for any $t \in \rz$. But this immediately implies the constancy of $\tilde{u}'$ and our claim follows. \hfill \qed\\

\begin{remark}\label{bern rem 1}
Let us end with a short remark on the failure of the Nitsche criterion (condition (4) in \cite{NN:1959_1})
for densities $g$ with the properties \reff{intro 7} and \reff{intro 8} from Theorem \ref{intro theo 2}. The boundedness of $g'$ in particular shows
that
\begin{equation}\label{bern 5}
\int_1^\infty g''(t) \dt < \infty \, , \quad \mbox{thus}\quad \int_1^\infty \frac{1}{\sqrt{s}}g''(\sqrt{s})\ds < \infty \, .
\end{equation}
Let us introduce the functions (compare \cite{NN:1959_1})
\[
f(t) := g (\sqrt{t}) \, , \quad \lambda(t) := \frac{2 f''(t)}{f'(t)} \, .
\]
Elementary calculations show
\begin{eqnarray*}
\frac{1 + t \lambda(t)}{2+t \lambda(t)} \cdot \frac{1}{t} &=& \frac{g''(\sqrt{t})}{\frac{1}{\sqrt{t}} g'(\sqrt{t}) + g''(\sqrt{t})} \cdot \frac{1}{t}
= \frac{1}{t + \frac{\sqrt{t}}{g''( \sqrt{t}) } g'(\sqrt{t})}
\end{eqnarray*}
and for $t \geq 1$ we obtain by the monotonicity of $g'$
\[
\frac{1+ t \lambda(t)}{2+t \lambda(t)} \frac{1}{t} \leq \frac{1}{g'(1)} \frac{g''(\sqrt{t})}{\sqrt{t}} \, .
\]
Recalling \reff{bern 5} it follows
\[
\int_1^\infty  \frac{1+ t  \lambda(t)}{2+t\lambda(t)} \frac{\dt}{t} < \infty \, ,
\]
which means that the ``Satz'' on p.~295 of \cite{NN:1959_1} does not apply to the linear growth case.

\end{remark}

%***************************************************************************************************
%***************************************************************************************************

\bibliography{Liouville_linear_extended}
\bibliographystyle{plain}

%\vspace*{0.5cm}
\begin{tabular}{lll}
Michael Bildhauer&e-mail:&bibi@math.uni-sb.de\\
Martin Fuchs&e-mail:&fuchs@math.uni-sb.de\\
Department of Mathematics&\\
Saarland University&\\
P.O.~Box 15 11 50&\\
66041 Saarbr\"ucken&\\ 
Germany
\end{tabular}
\end{document}